\documentclass[11pt, a4paper]{amsart}

\usepackage{amssymb,amsmath,amsthm}
\usepackage{hyperref}
\usepackage{graphicx}
\usepackage{xcolor}
\usepackage{graphicx} 
\hypersetup{backref, pdfpagemode=FullScreen, colorlinks=true,
  citecolor=magenta, linkcolor=cyan, urlcolor=blue}

\usepackage{mathtools}
\mathtoolsset{showonlyrefs}

\newcommand{\Kon}{\mathcal{K}_o^n}
\newcommand{\Ken}{\mathcal{K}_e^n}
\newcommand{\Kn}{\mathcal{K}^n}

\newcommand{\vol}{\mathrm{vol}}

\newcommand{\relint}{\mathrm{relint}\,}
\newcommand{\conv}{\mathrm{conv}}

\newcommand{\lin}{\mathrm{lin}\,}
\newcommand{\R}{\mathbb{R}}

\newcommand{\N}{\mathbb{N}}

\newcommand{\W}{\mathrm{W}}
\newcommand{\cura}{\mathrm{C}}
\newcommand{\sura}{\mathrm{S}}
\newcommand{\dW}{\widetilde{\mathrm{W}}}
\newcommand{\dcura}{\widetilde{\mathrm{C}}}

\newcommand{\va}{{\boldsymbol a}}
\newcommand{\vb}{{\boldsymbol b}}

\newcommand{\ve}{{\boldsymbol e}}

\newcommand{\vt}{{\boldsymbol t}}

\newcommand{\vu}{{\boldsymbol u}}
\newcommand{\vv}{{\boldsymbol v}}

\newcommand{\vx}{{\boldsymbol x}}
\newcommand{\vy}{{\boldsymbol y}}
\newcommand{\vz}{{\boldsymbol z}}

\newcommand{\vnull}{{\boldsymbol 0}}
\newcommand{\ov}{\overline}
\DeclareMathOperator{\dplus}{{\widetilde +}}
\newtheorem{theorem}{Theorem}[section]
\newtheorem{corollary}[theorem]{Corollary}
\newtheorem{lemma}[theorem]{Lemma}
\newtheorem{remark}[theorem]{Remark}

\newtheorem{proposition}[theorem]{Proposition}

\numberwithin{equation}{section}

\usepackage{pdfcomment}

\begin{document}

\title[Subspace concentration  of  dual curvature measures]{Subspace
  concentration  of  dual curvature measures of symmetric
  convex bodies}
\author{K\'aroly J. B\"or\"oczky}
\address{Alfr\'ed R\'enyi Institute of Mathematics,
		Hungarian Academy of Sciences,
		Re\'altanoda u. 13-15.,
		H-1053 Budapest,
		Hungary}
\email{carlos@renyi.hu}
\author{Martin Henk}
\author{Hannes Pollehn}
\address{Technische Universit\"at Berlin,
		Institut f\"ur Mathematik,
		Strasse des 17. Juni 136,
		D-10623 Berlin,
		Germany}
\email{henk@math.tu-berlin.de, pollehn@math.tu-berlin.de}
%\thanks{}
%\dedicatory{}
\begin{abstract}
We prove a tight subspace concentration inequality 
for the dual curvature measures of a symmetric convex body. 
\end{abstract}
\maketitle

\section{Introduction}
Let $\Kn$ denote the set of convex bodies in $\R^n$,
i.e., all convex and compact subsets $K$ having a non-empty
interior. The set of convex bodies having the origin as an interior
point and the set of origin-symmetric convex bodies, i.e., those sets
which satisfy $K=-K$ 
are denoted by $\Kon$ and $\Ken$ respectively. For $\vx,\vy\in\R^n$,
let $\langle \vx,\vy\rangle$ denote the standard inner product and
$|\vx|=\sqrt{\langle \vx,\vx\rangle}$ the Euclidean norm. We write
$B_n$ for the $n$-dimensional Euclidean unit ball, i.e., 
$B_n=\{ \vx\in\R^n : |\vx|\leq 1 \}$ and $S^{n-1}$ for its
boundary.  The $k$-dimensional Hausdorff-measure will be denoted by
$\mathcal{H}^{k}(\cdot)$ and instead of $\mathcal{H}^{n}(\cdot)$ we
will also write $\vol(\cdot)$ for the $n$-dimensional volume. 

At the heart of the Brunn-Minkowski theory is the study of the volume
functional with respect to the Minkowski addition of convex bodies. This leads to the theory
of mixed volumes and, in particular, to  the quermassintegrals $\W_i(K)$ of a
convex body $K\in\Kn$.  The latter  may be defined via the classical Steiner
formula, expressing the volume of the Minkowski sum of $K$ and
$\lambda\,B_n$, i.e., the volume of the parallel body of $K$ at distance
$\lambda$ 
%  outer parallel body at distance
% $\lambda$, i.e.  
% \begin{equation}
%    K+\lambda\,B_n=\{\vx\in\R^n : \mathrm{d}(K,\vx)\leq \lambda\}
% \end{equation} 
as a polynomial in $\lambda$  (cf., e.g., \cite[Sect.~4.2]{Schneider:1993})
\begin{equation}
 \vol(K+\lambda\,B_n)=\sum_{i=0}^n \lambda^i\, \binom{n}{i} \W_i(K).
\label{eq:steiner} 
\end{equation}  
A more direct geometric interpretation is given
by Kubota's integral formula (cf., e.g., \cite[Subsect.~5.3.2]{Schneider:1993}),  showing that they are -- up to some
constants  -- the means of the volumes of projections 
\begin{equation} 
    \W_{n-i}(K)=\frac{\vol(B_n)}{\vol_i(B_i)}\int_{G(n,i)}\vol_i(K|L)\,\mathrm{d}L,
    \quad i=1,\dots,n, 
\label{eq:kubota}
\end{equation}  
where $\vol_i(\cdot)$ denotes the $i$-dimensional volume, integration
is taken with respect to the rotation-invariant probability measure on
the Grassmannian $G(n,i)$ of all $i$-dimensional linear subspaces  and
$K|L$ denotes the orthogonal projection onto $L$.  

A local version of the Steiner formula above leads to two important
series of geometric measures, the  
surface area measures $\sura_i(K,\cdot)$ and the curvature measures
$\cura_i(K,\cdot)$, $i=0,\dots,n-1$, of a convex body $K$. Here we will
only briefly describe the surface area measures since they may be
considered as the ``primal''  counterpart to  the dual curvature measures we
are interested in.  

To this end, we denote for 
$\omega\subseteq S^{n-1}$ by $\nu^{-1}_K(\omega)\subseteq\partial K$  the set of all boundary points of $K$ having
an outer unit normal in $\omega$.   Moreover, for
$\vx\in\R^n\setminus K$ let $r_K(\vx)\in \partial K$ be the point of $K$
closest to $K$.  Then for a Borel set
$\omega\subseteq S^{n-1}$ and $\lambda>0$ we consider the 
local  parallel body 
\begin{equation}
  B_K(\lambda,\omega) = \left\{\vx\in\R^n : 0< |\vx-r_K(\vx)|\leq \lambda
  \text{ and } r_K(\vx)\in\nu_K^{-1}(\omega)\right\}. 
\label{eq:local_parallel}
\end{equation}  
The local Steiner formula is now a polynomial in $\lambda$ 
 whose  coefficients are (up to constants
 depending on $i,n$)  the surface area measures (cf., e.g., \cite[Sect.~4.2]{Schneider:1993})
\begin{equation}
\vol( B_K(\lambda,\omega)) = \frac{1}{n} \sum_{i=1}^{n} \lambda^i\,
\binom{n}{i} \sura_{n-i}(K,\omega).
\label{eq:local_steiner} 
\end{equation}
They may also be regarded as the (right hand side) differentials of
the quermassintegrals 
\begin{equation}
 \lim_{\epsilon \downarrow 0}\frac{\W_{n-1-i}(K+\epsilon
   B_n)-\W_{n-1-i}(K)}{\epsilon}=\int_{S^{n-1}}\,\mathrm{d}\sura_{i}(K,\vu).
\label{eq:differentials}
\end{equation} 
Also observe that $\sura_i(K,S^{n-1})=n\,\W_{n-i}(K)$, $i=0,\dots, n-1$.  

To characterize the surface area measures $\sura_i(K,\cdot)$,
 $i\in\{1,\dots,n-1\}$,  among the  finite Borel
measures on the sphere is a corner stone of
the Brunn-Minkowski theory. Today this  problem is known as  the {\em Minkowski--Christoffel} 
problem, since for $j=n-1$ and the surface area measure
$\sura_{n-1}(K,\cdot)$ it is the classical Minkowski problem and for $j=1$
it is the Christoffel problem.  We refer to \cite[Chapter
8]{Schneider:1993} for more information and references.

There are two far-reaching extensions of the classical Brunn-Minkowski theory,
both arising 
basically  by replacing the classical Minkowski-addition
by another additive operation (cf.~\cite{Gardner:2014}). The first one is the $L_p$ addition
introduced by Firey (see, e.g.,  \cite{Firey:1962}) which leads to the rich and emerging 
{\em $L_p$-Brunn-Minkowski theory} for which we refer to
\cite[Sect.~9.1, 9.2]{Schneider:1993}). 

The second one, introduced by Lutwak \cite{Lutwak:1975a},  is based on the radial addition $\dplus$ where $\vx
\dplus\vy = \vx+\vy$  if $\vx,\vy$ are linearly dependent and
$\vnull$ otherwise.  Considering the volume of radial additions leads
to the {\em dual Brunn-Minkowski theory} (cf.~\cite[Sect.~9.3]{Schneider:1993}) 
with dual mixed volumes, and, in
particular, also with dual quermassintegrals $\dW_i(K)$
arising via a dual Steiner formula  (cf.~\eqref{eq:steiner})
\begin{equation}
 \vol(K\dplus\lambda\,B_n)=\sum_{i=0}^n \lambda^i\, \binom{n}{i}
 \dW_i(K).
\label{eq:dual_steiner} 
\end{equation} 
In general the radial addition of two convex sets is not a convex set,
but the radial addition of two star bodies is again a  star body. This is one of the features of the dual
Brunn-Minkowski theory which makes it so useful.    
The celebrated solution
of the Busemann-Petty problem is amongst the recent successes 
of the  dual Brunn-Minkowski theory, cf. \cite{Gardner:1994,
  GardnerKoldobskySchlumprecht:1999, Zhang:1999}, and it 
 also has connections and applications to integral geometry,
Minkowski geometry, and the local theory of Banach spaces. 

In analogy to Kubota's formula \eqref{eq:kubota}  the dual quermassintegrals $\dW_i(K)$  admit the following integral
geometric representation as the means of the volumes of sections (cf.~\cite[Sect. 9.3]{Schneider:1993}) 
\begin{equation}
    \dW_{n-i}(K)=\frac{\vol(B_n)}{\vol_i(B_i)}\int_{G(n,i)}\vol_i(K\cap
    L)\,\mathrm{d}L,
    \quad i=1,\dots,n. 
\end{equation} 
There are many more ``dualities'' between the classical and
dual theory, but  there were no dual geometric
measures corresponding to the  surface area or curvature measures.
This missing link was recently established in the 
ground-breaking paper \cite{Huang:2016} by Huang, Lutwak, Yang and
Zhang. Let $\rho_K$ be the radial function (see Section 2 for the definition)
  of a convex body $K\in \Kon$. % For $\vx \in\R^n\setminus K$
  % let $\widetilde{d}(K,\vx)$ be the Euclidean distance between $\vx$
  % and the boundary point $\rho_K(\vx)\vx\in\partial K$.
  Analogous to
  \eqref{eq:local_parallel} we consider for a Borel set $\eta\subseteq S^{n-1}$ and $\lambda>0$ the set 
\begin{equation}
 {\widetilde A}_K(\lambda,\eta) = \left\{\vx\in\R^n : 0\leq|\vx-\rho_K(\vx)\vx|\leq \lambda
  \text{ and } \rho_K(\vx)\vx\in\nu_K^{-1}(\eta)\right\}. 
\end{equation}    
Then there also exists  a local Steiner type formula of these  local
dual parallel sets \cite[Theorem 3.1]{Huang:2016} (cf.~\eqref{eq:local_steiner})
\begin{equation}
        \vol({\widetilde A}_K(\lambda,\eta))=\sum_{i=0}^n \binom{n}{i}\lambda^{i}\dcura_{n-i}(K,\eta).
\end{equation} 
$\dcura_{i}(K,\eta)$ is called the {\em $i$th dual curvature
  measure} and they are the counterparts to the surface area measures
$\sura_i(K,\omega)$  in the dual Brunn-Minkowski theory. Observe that
$\dcura_i(K,S^{n-1})=\dW_{n-i}(K)$. 
As the surface area measure (cf.~\eqref{eq:differentials}), the dual
curvature measures  may also be
considered as differentials of the dual quermassintegrals, even in a
stronger form (see  \cite[Section 4]{Huang:2016}). 
We want to 
point out  that  there are also
dual surface area measures corresponding to the curvature measures
in the classical theory (see \cite{Huang:2016}).   

Huang, Lutwak, Yang and Zhang also gave  an explicit 
integral representation of the dual curvature measures  which allowed
them to define more generally for $q\in \R$ the  $q$th dual
curvature measure  of a convex body $K\in\Kon$ as \cite[Def. 3.2]{Huang:2016} 
\begin{equation}
	\dcura_q(K,\eta) = \frac1n
        \int\limits_{\alpha_K^\ast(\eta)} \rho_K(\vu)^q
        \mathrm{d}\mathcal{H}^{n-1}(\vu). 
\label{eq:dual_curvature_measure}
\end{equation}
Here  $\alpha_K^\ast(\eta)$
denotes the set of directions $\vu\in S^{n-1}$, such that the boundary
point $\rho_K(\vu)\vu$ belongs to $\nu_K^{-1}(\eta)$. The analog to the Minkowski-Christoffel problem in the dual
Brunn-Minkowski theory is  (cf.~\cite[Sect.~5]{Huang:2016}) 

\smallskip
\begin{center}
\begin{minipage}{0.95\textwidth}{\em The  dual Minkowski problem.} Given  a finite Borel measure $\mu$ on
  $S^{n-1}$ and $q\in \R$. Find necessary and sufficient conditions
  for the existence of a convex body $K\in\Kon$ such that $\dcura_q(K,\cdot)=\mu$.  
\end{minipage} 
\end{center}
\smallskip 
%  These
% measures  arise surprisingly naturally in the dual
% Brunn-Minkowski theory and may be considered as the dual
% counterparts to Federer's area measures.  In particular, for
% $q\in\{0,\dots,n\}$ we have
% \begin{equation} 
% \dcura_q(K,S^{n-1})=\widetilde\W_{n-q}(K),  
% \end{equation}  
% where $\widetilde\W_{n-q}(K)$ is the $(n-q)$-th dual
% quermassintegral (see \cite{Lutwak:1975a, Lutwak:1975b})

\noindent
An amazing feature of these dual curvature measures
is that they also link two other well-known  fundamental geometric measures of a convex body
(cf.~\cite[Lemma 3.8]{Huang:2016}):
when  $q=0$  the dual curvature measure
$\dcura_0(K,\cdot)$ is -- up to a factor of $n$ -- 
Aleksandrov's integral curvature of the polar body of $K$ and for
$q=n$ the dual curvature measure coincides with the cone-volume
measure of $K$ given by 
\begin{equation}
\dcura_n(K,\eta)	= V_K(\eta) = \frac{1}{n} \int\limits_{\nu_K^{-1}(\eta)} \langle \vx,\nu_K(\vx)\rangle \mathrm{d}\mathcal{H}^{n-1}(\vx).
\end{equation}

% where  $\nu_K:\partial'K\to S^{n-1}$ is the Gauss map, i.e.,
% $\partial'K$, is the set of all points $\vx$ in the boundary of $K$
% having a unique outer normal vector $\nu_K(\vx)\in S^{n-1}$. 

%%% TEXT BEGIN

Similarly to the Minkowsi problem, the dual Minkowski problem is equivalent to solving a Monge-Amp\`{e}re type partial differential equation if the measure $\mu$ has a density function $g : S^{n-1}\to \R$. In particular, if $q\in(0,n]$, then
the dual Minkowski problem amounts to solving the Monge-Amp\`ere  equation
\begin{equation}
\label{Monge-Ampere}
\frac1n\,h(\vx)|\nabla h(\vx) + h(\vx)\,\vx |^{q-n}\det[h_{ij}(\vx) + \delta_{ij}h(\vx)] = g(\vx),
\end{equation}
where $[h_{ij}(\vx)]$ is the Hessian matrix of the (unknown) support function $h$ with respect to an orthonormal
frame on $S^{n-1}$, and $\delta_{ij}$ is the Kronecker delta.

If $\frac1n\,h(\vx)|\nabla h(\vx) + h(\vx)\vx|^{q-n}$ were omitted in \eqref{Monge-Ampere}, then \eqref{Monge-Ampere} would become the partial differential
equation of the classical Minkowski problem, see, e.g., 
\cite{Caffarelli:1990,Cheng:1976,Nirenberg:1953}. If only the factor
$|\nabla h(\vx) + h(\vx)\,\vx|^{q-n}$ were omitted, then
equation \eqref{Monge-Ampere} would become the partial differential equation associated with the 
cone volume measure, the so-called logarithmic Minkowski problem (see,
e.g., \cite{Boroczky:2013, Chou:2006}). Due to the gradient component
in \eqref{Monge-Ampere} if $q\in(0,n)$, the dual Minkowski problem is signicantly more challenging than the classical Minkowski problem and logarithmic Minkowski problem.

%%% TEXT END

The cone-volume measure for convex bodies  has been studied extensively over the last few
years in many different contexts, see, e.g., \cite{Barthe:2005,
  Boroczky:2015c, Boroczky:2012,Boroczky:2013,Boroczky:2015a,
  Gardner:2014,  Gromov:1987, Haberl:2014, Huang:2016, Ludwig:2010a, Ludwig:2010b,
  Lutwak:2005, Lutwak:2010a, Lutwak:2010b, Ma:2015, Naor:2007, Naor:2003, Paouris:2012,
  Stancu:2012, Zhu:2014, Zhu:2015}.  
% ; see for instance 
% Barthe, Guedon, Mendelson and  Naor~\cite{Barthe:2005}, 
% B\"or\"oczky and P. Heged{\H u}s~\cite{Boroczky:2015c},
% B\"or\"oczky, Lutwak, Yang and Zhang~\cite{Boroczky:2012,Boroczky:2013,Boroczky:2015a}, 
% Gardner, Hug and Weil~\cite{Gardner:2014},
% Gromov and Milman~\cite{Gromov:1987}, 
% Huang, Lutwak, Yang and Zhang \cite{Huang:2016},
% Ludwig~\cite{Ludwig:2010a}, 
% Ludwig and Reitzner~\cite{Ludwig:2010b}, 
% Lutwak, Yang and Zhang~\cite{Lutwak:2005}, 
% Ma~\cite{Ma:2015},
% Naor~\cite{Naor:2007}, 
% Naor and Romik~\cite{Naor:2003}, 
% Paouris and Werner~\cite{Paouris:2012}, 
% Stancu~\cite{Stancu:2012}, 
% Zhu~\cite{Zhu:2014, Zhu:2015}.
One very important property of the cone-volume measure -- and which
makes it is so useful --  is its
$\mathrm{SL}(n)$-invariance, or simply called affine invariance.  It
is also the subject of the central {\em logarithmic Minkowski problem}
which asks  for sufficient and necessary conditions 
of a measure $\mu$ on $S^{n-1}$  to be  the cone-volume measure of a
convex body $K\in\Kon$. This is the $p=0$ limit case of the general
$L_p$-Minkowski problem within the above mentioned $L_p$
Brunn-Minkowski theory for which we refer to \cite{Hug:2005,
  Lutwak:1993, Zhu:2015a} and the references within.

The discrete, planar, even case of the logarithmic Minkowski
problem, i.e., with respect to origin-symmetric convex polygons, was
completely solved by Stancu~\cite{Stancu:2002,Stancu:2003}, and later
Zhu~\cite{Zhu:2014} as well as B\"or\"oczky, Heged{\H{u}}s and
Zhu~\cite{Boroczky:2015b} settled (in particular) the  case when $K$ is a polytope 
whose outer normals are in general position.

 % In~\cite{Boroczky:2013} it was shown, that the following condition is necessary and sufficient for even measures to be a solution to the logarithmic Minkowski problem when $K$ is symmetric.
In~\cite{Boroczky:2013}, B\"or\"oczky, Lutwak, Yang and Zhang gave a
complete characterization of the cone-volume measure of
origin-symmetric convex bodies among the even measures on the sphere. 
The key feature of such a measure is expressed via the following
condition:  A non-zero,
finite Borel measure $\mu$ on the unit sphere satisfies the 
\emph{subspace concentration condition} if
\begin{equation}
	\frac{\mu(S^{n-1}\cap L)}{\mu(S^{n-1})} \leq \frac{\dim L}{n}
	\label{eq:subspace_concentration_inequality}
\end{equation}
for every proper subspace $L$ of $\R^n$, and whenever we have equality in~\eqref{eq:subspace_concentration_inequality} for some $L$, there is a subspace $L'$ complementary to $L$, such that $\mu$ is concentrated on $S^{n-1}\cap (L\cup L')$.

Apart from the uniqueness aspect of the Minkowski problem, 
the symmetric case of the logarithmic Minkowski problem is settled.
\begin{theorem}[\cite{Boroczky:2013}]
	A non-zero, finite, even Borel measure $\mu$ on $S^{n-1}$ is
        the cone-volume measure of $K\in \Ken$  if and only if $\mu$ satisfies the subspace concentration condition.
	\label{thm:logarithmic_minkowski_problem}
\end{theorem}
% The necessity of the inequality \eqref{eq:subspace_concentration_inequality} for
% cone-volume measures of symmetric polytopes was earlier independently
% proven by Henk, Sch\"urmann and Wills~\cite{Henk:2005} (in a different
% context) and by He, Leng and Li~\cite{He:2006}. 
An extension of the
validity of  inequality \eqref{eq:subspace_concentration_inequality}
to centered bodies, i.e.,   bodies whose center of mass is at the
origin, was given  in the discrete case by Henk and Linke \cite{Henk:2014}, and
in the general setting by B\"or\"oczky and Henk~\cite{Boroczky:2016}.

% In their very recent and ground-breaking paper Huang, Lutwak, Yang and Zhang~\cite{Huang:2016} introduced the dual curvature measures of convex bodies, which arise naturally in the dual Brunn-Minkowski theory. For $K\in\Kon$ and $\eta\subseteq S^{n-1}$ let $\alpha_K^\ast(\eta)$ be the set of directions $u\in S^{n-1}$, such that the boundary point $\rho_K(u)u$ has an outer normal lying in $\eta$. For a real number $q\in\R$ we define the $q$-th dual curvature measure $\dcura_q(K,\cdot)$ of $K$ by
% \begin{equation}
% 	\dcura_q(K,\eta) = \frac1n \int\limits_{\alpha_K^\ast(\eta)} \rho_K(u)^q \mathrm{d}\mathcal{H}^{n-1}(u)
% 	\label{eq:dual_curvature_measure}
% \end{equation}
% for every Borel set $\eta\subseteq S^{n-1}$. For $q=n$ the curvature
% measure coincides with the cone-volume measure.

A generalization (up to the equality case) of the sufficiency part of Theorem
\ref{thm:logarithmic_minkowski_problem} to the $q$-dual curvature
measure  for
$q\in(0,n]$ was  given by Huang, Lutwak, Yang and Zhang.
For clarity, we separate their main result into 
the next two theorems.

\begin{theorem}[\protect{\cite[Theorem 6.6]{Huang:2016}}]
If $q\in (0,1]$, then an even finite Borel measure $\mu$ on
$S^{n-1}$ is a $q$-dual curvature measure if and only if $\mu$ is not concentrated on any great subsphere.
\label{thm:subspace_mass0}
\end{theorem}

\begin{theorem}[\protect{\cite[Theorem 6.6]{Huang:2016}}]
Let $q\in [1,n]$ and let  $\mu$  be a non-zero, finite,  even Borel measure  on
$S^{n-1}$ satisfying  the subsapce mass inequality 
\begin{equation}
	\frac{\mu(S^{n-1}\cap L)}{\mu(S^{n-1})} <
        1-\frac{q-1}{q}\frac{n-\dim L}{n-1}
\label{eq:subspace_mass} 
\end{equation}
for every proper subspace $L$ of $\R^n$. Then there exists a
$o$-symmetric convex body  $K\in\Ken$ with $\dcura_q(K,\cdot)=\mu$. 
\label{thm:subspace_mass}
\end{theorem}

In particular, it is highly desirable to understand how close 
\eqref{eq:subspace_mass} is to characterize $q$-dual curvature measures.
Observe that for $q=n$ the inequality \eqref{eq:subspace_mass} becomes
essentially \eqref{eq:subspace_concentration_inequality}.

Our main result treats the necessity of a subspace concentration bound on
dual curvature measures.  % of symmetric convex bodies similar to the ones stated above. We also provide examples showing, that the bounds in the following theorem cannot be lowered.
\begin{theorem}
Let $K\in\Ken$, $q\in [1,n]$ and let $L\subset \R^n$ be a proper
subspace. Then we have 
\begin{equation}
	\frac{\dcura_q(K,S^{n-1}\cap L)}{\dcura_q(K,S^{n-1})} \leq
        \min\left\{ \frac{\dim L}{q},1\right\}, 
	\label{eq:subspace_bound}
\end{equation}
%for every proper linear subspace $L\subset\R^n$. 
and equality holds in \eqref{eq:subspace_bound} if and only if $q=n$
and $\dcura_n(K,\cdot)$, i.e., the cone-volume measure of $K$,  satisfies the subspace
concentration condition \eqref{eq:subspace_concentration_inequality}.
\label{thm:subspace_bound}
\end{theorem}

In particular, for $q<n$ we always have strict inequality in
\eqref{eq:subspace_bound}, but this is also optimal. 
\begin{proposition} Let $0<q<n$ and $k\in\{1,\dots,n-1\}$.  There exists a sequence of convex
  bodies  $K_l\in\Ken$, $l\in\N$, and a $k$-dimensional subspace
  $L\subset\R^n$ such that 
\begin{equation}
       \lim_{l\to\infty} \frac{\dcura_q(K_l,S^{n-1}\cap
         L)}{\dcura_q(K_l,S^{n-1})}=\begin{cases} 
  \frac{k}{q}&, k\leq q, \\ 
1&,k\geq q.
\end{cases} 
\end{equation} 
\label{prop:bound_opt}  
\end{proposition}

We observe that if $q\in[1,n]$ and $\dim L=1$ for a linear subspace $L$ then  
$1-\frac{q-1}{q}\frac{n-\dim L}{n-1}=\frac{\dim L}{q}$. Therefore 
Theorems~\ref{thm:subspace_mass} and \ref{thm:subspace_bound} 
complete the characterization of the $q$-dual curvature measures if $n=2$.

\begin{corollary}
If $q\in [1,2)$, then an even finite Borel measure $\mu$ on
$S^1$ is a $q$-dual curvature measure if and only if 
$$
\frac{\mu(S^1\cap L)}{\mu(S^1)} <\frac{1}{q}
$$
for every one-dimensional subspace $L$ of $\R^2$.
\label{cor:subspace_mass-planar}
\end{corollary}

We remark that the dual
Minkowski problem is far easier to handle for the special case 
where the measure $\mu$  has a positive continuous density, (where
subspace concentration is trivially satisfied). The singular general
case for measures is substantially more delicate, which involves
measure concentration and requires far more powerful techniques to
solve.

The paper is organized as follows. First we will briefly recall
some basic facts about convex bodies needed in our investigations in
Section 2.  In
Section 3 we will prove a lemma in the spirit of the celebrated
Brunn-Minkowski theorem, which is one of the main ingredients for the
proof of Theorem~\ref{thm:subspace_bound} given in Section 4. Finally,
in Section 5 we will prove Proposition \ref{prop:bound_opt}. 

\section{Preliminaries}
We recommend the books by Gardner \cite{Gardner:2006}, Gruber \cite{Gruber:2007} and
Schneider \cite{Schneider:1993} as excellent references on  convex geometry.

For a given convex body $K\in\Kn$ the support function $h_K\colon\R^n\to\R$ is defined by
\begin{equation}
	h_K(\vx)=\max_{\vy\in K} \langle \vx,\vy\rangle.
\end{equation}
A boundary point $\vx\in\partial K$ is said to have a (not necessarily unique) unit outer normal vector $\vu\in S^{n-1}$ if $\langle \vx,\vu\rangle = h_K(\vu)$. The corresponding supporting hyperplane $\{\vx\in \R^n\colon \langle \vx,\vu\rangle=h_K(\vu)\}$ will be denoted by $H_K(\vu)$. For $K\in\Kon$ the radial function $\rho_K\colon \R^n\setminus\{\vnull\}\to\R$ is given by
\begin{equation}
	\rho_K(\vx)=\max\{\rho>0\colon \rho\, \vx\in K\}.
\end{equation}
Note, that the support function and the radial function are
homogeneous of degrees $1$ and $-1$, respectively, i.e., 
\begin{equation} 
   h_K(\lambda\,\vx) =\lambda\, h_K(\vx)\text{ and } \rho_K(\lambda\,\vx)=\lambda^{-1}\, \rho_K(\vx),
\end{equation} 
for $\lambda>0$.
  We define the reverse radial Gauss image of $\eta\subseteq S^{n-1}$
  with respect to a convex body $K\in\Kon$ by
\begin{equation}
	\alpha_K^\ast(\eta) = \{\vu\in S^{n-1}\colon \rho_K(\vu)\vu\in H_K(\vv) \text{ for a }\vv\in\eta \}.
\end{equation}
If $\eta$ is a Borel set, then $\alpha_K^\ast(\eta)$ is
$\mathcal{H}^{n-1}$-measurable
(see~\cite[Lemma~2.2.11.]{Schneider:1993}) and so the  $q$-th dual
curvature measure given in~\eqref{eq:dual_curvature_measure} is well defined. We will need the following identity.
\begin{lemma}
	Let $K\in\Kon$, $q>0$ and $\eta\subseteq S^{n-1}$ a Borel set.
        Then 
	\begin{equation}
	\dcura_q(K,\eta) = \frac{q}{n} \int\limits_{\vx\in K,\,\,\vx/|\vx|\in\alpha_K^\ast(\eta)} |\vx|^{q-n} \mathrm{d}\mathcal{H}^n(\vx).
	\end{equation}
	\label{lem:curvature_measure_euclidean_coordinates}
\end{lemma}
\begin{proof}
	Since $q>0$, we may write via spherical coordinates 
\begin{equation}
\begin{split}
 \frac{q}{n} \int\limits_{\vx\in
   K,\,\,\vx/|\vx|\in\alpha_K^\ast(\eta)} |\vx|^{q-n}
 \mathrm{d}\mathcal{H}^n(\vx) & = 
 \frac{1}{n}\int\limits_{\alpha_K^\ast(\eta)}\left(\int\limits_0^{\rho_K(\vu)}
   q\,r^{n-1\,}r^{q-n} \mathrm{d}r\right)
 \mathrm{d}\mathcal{H}^{n-1}(\vu)\\
&= \frac1n
        \int\limits_{\alpha_K^\ast(\eta)} \rho_K(\vu)^q
        \mathrm{d}\mathcal{H}^{n-1}(\vu)= 	\dcura_q(K,\eta).
\end{split}  
\end{equation} 
% 	\begin{equation}
% \begin{split} 
% 	\dcura_q(K,\eta) &=\frac1n
%         \int\limits_{\alpha_K^\ast(\eta)} \rho_K(\vu)^q
%         \mathrm{d}\mathcal{H}^{n-1}(\vu) \\ &= \frac{1}{n}
%         \int\limits_{\alpha_K^\ast(\eta)}
%         \left(\int\limits_0^{\rho_K(\vu)} q\,r^{q-1} \mathrm{d}r\right)
%         \mathrm{d}\mathcal{H}^{n-1}(\vu).
% \end{split}
% 	\end{equation}
% 	which by going from spherical to Euclidean coordinates via $x=r\cdot u$, $\mathrm{d}\mathcal{H}^n(x)=r^{n-1}\mathrm{d}r\mathrm{d}\mathcal{H}^{n-1}(u)$, becomes
% 	\begin{equation}
% 	\dcura_q(K,\eta) = \frac{1}{n} \int\limits_{\substack{x\in K\\ x/|x|\in\alpha_K^\ast(\eta)}} q|x|^{q-1}\cdot |x|^{1-n} \mathrm{d}\mathcal{H}^n(x).
% 	\end{equation}
\end{proof}

Let $L$ be a linear subspace of $\R^n$. We write $K|L$ to denote the
orthogonal projection of $K$ onto $L$  and $L^\perp$ for the subspace orthogonal to $L$.

As usual, for two subsets  $A, B\subseteq \R^n$ and reals
$\alpha,\beta \geq 0$ the  Minkowski combination is defined by 
\begin{equation}
	\alpha A+\beta B = \{\alpha \va+\beta \vb : \va\in A,\vb\in
        B\}. 
\end{equation} 
% and we also write $\alpha A-\beta B$ instead of  $\alpha A + (-\beta)
% B = \alpha A + \beta (-B)$.

By the well-known Brunn-Minkowski inequality we know that the $n$-th root 
of the volume of the Minkowski combination  is a concave function. 
More precisely, for two convex bodies
$K_0,K_1\subset\R^n$ and for $\lambda\in[0,1]$  we have 
\begin{equation}
	\vol_n((1-\lambda)K_0+\lambda K_1)^{1/n}\geq (1-\lambda)\vol_n(K_0)^{1/n}+ \lambda \vol_n(K_1)^{1/n},
	\label{eq:brunn_minkowski}
\end{equation}
where $\vol_n(\cdot)=\mathcal{H}^n(\cdot)$ denotes the $n$-dimensional
Hausdorff measure.  We have equality in \eqref{eq:brunn_minkowski}
for some $0<\lambda<1$ if and only if $K_0$ and $K_1$ lie in parallel
hyperplanes  or they are homothetic, i.e., there exist a $\vt\in\R^n$ and $\mu\geq 0$ such that $K_1=\vt+\mu\,K_0$ (see, e.g., \cite{Gardner:2002}, \cite[Sect. 6.1]{Schneider:1993}).

\section{Integrals of even unimodal functions}
A function $f$ on the real line $\R$ is called unimodal if there is a
number $m\in\R$, 
such that $f$ is an increasing function on $(-\infty,m)$ 
and decreasing on $(m,\infty)$. Obviously the integral of $f$ over an
interval of fixed length is maximal when the interval is centered at
$m$. The notion of unimodal functions can be extended to 
higher dimensional spaces in the following way. 

The superlevel sets of a function $f:\R^n\to\R$ are given by
$L_f^+(\alpha)=\{ \vx\in\R^n\colon f(\vx)\geq\alpha \}$,
$\alpha\in\R$. We say that $f$ is unimodal if every superlevel set of
$f$ is closed and convex. 
It was shown by Anderson \cite{Anderson:1955} that the integral of an
even unimodal function  over translates of a symmetric convex region 
is maximal if the center of symmetry is moved to the origin. 
His proof relies only on the Brunn-Minkowski theorem. Here we
generalize  this approach to integrals over a convex combination of a convex body $K$ and its reflection $-K$.
\begin{lemma}
Let $K\in\Kn$ and $\dim K=k$. Let $f\colon\R^n\to\R_{\geq
  0}\cup\{\infty\}$ be a unimodal function,  such that $f(\vx)=f(-\vx)$ for every $\vx\in\R^n$ and
\begin{equation}
	\int\limits_{\frac12 K+\frac12(-K)} f(\vx) \mathrm{d}\mathcal{H}^k(\vx)<\infty.
\end{equation}
Let $\lambda\in (0,1)$. Then
\begin{equation}
	\int\limits_{\lambda K+(1-\lambda)(-K)} f(\vx)\mathrm{d}\mathcal{H}^k(\vx) \geq \int\limits_K f(\vx) \mathrm{d}\mathcal{H}^k(\vx).
	\label{eq:integral_conv_comb}
\end{equation}
Moreover, equality holds if and only if for every $\alpha>0$
\begin{equation}
	\vol_k\left([\lambda K+(1-\lambda)(-K)]\cap L_f^+(\alpha)\right) = \vol_k(K\cap L_f^+(\alpha)).
\end{equation}
\label{lem:unimodal}
\end{lemma}
\begin{proof}
Let $K_\lambda = \lambda K+(1-\lambda)(-K)$. By the convexity of
$L_f^+(\alpha)$ we have for every $\alpha\in\R$ 
\begin{equation}
	K_\lambda \cap L_f^+(\alpha)\supseteq \lambda (K\cap L_f^+(\alpha)) + (1-\lambda) ((-K)\cap L_f^+(\alpha)).
	\label{eq:level_sets_conv_comb}
\end{equation}
The Brunn-Minkowski inequality \eqref{eq:brunn_minkowski} applied to
the  set on right hand side of \eqref{eq:level_sets_conv_comb} gives 
\begin{equation}
\begin{split}
\vol_k(K_\lambda\cap L_f^+(\alpha)) \geq
&\vol_k\left( \lambda (K\cap L_f^+(\alpha)) + (1-\lambda) ((-K)\cap L_f^+(\alpha)) \right)\\
\geq& \left(\lambda \vol_k(K\cap L_f^+(\alpha))^{1/k}+(1-\lambda)\vol_k((-K)\cap L_f^+(\alpha))^{1/k}\right)^k.
\end{split}
\end{equation}
Since $f$ is even, the superlevel sets $L_f^+(\alpha)$ are symmetric. 
Hence,  $\vol_k(K\cap L_f^+(\alpha))=\vol_k((-K)\cap L_f^+(\alpha))$
and so 
\begin{equation}
	\vol_k(K_\lambda\cap L_f^+(\alpha)) \geq \vol_k(K\cap L_f^+(\alpha))
\end{equation}
 for every $\alpha\in\R$. Fubini's theorem yields
\begin{align*}
	\int\limits_{K_\lambda} f(x)\mathrm{d}\mathcal{H}^k(x) &= \int\limits_0^\infty \vol_k(K_\lambda\cap L_f^+(\alpha))\mathrm{d}\alpha\\
	&\geq \int\limits_0^\infty \vol_k(K\cap L_f^+(\alpha))\mathrm{d}\alpha\\
	&= \int\limits_K f(x) \mathrm{d}\mathcal{H}^k(x).
\end{align*}
Suppose we have equality in \eqref{eq:integral_conv_comb}. Since
$\vol_k(K_\lambda\cap L_f^+(\alpha))$ is  continuous on the left with
respect to $\alpha$ we find that
\begin{equation}
\vol_k(K_\lambda\cap L_f^+(\alpha)) = \vol_k(K\cap L_f^+(\alpha))
\end{equation}
for every $\alpha>0$.
\end{proof}

\section{Proof of Theorem~\ref{thm:subspace_bound}}
Now we are ready to give the proof of
Theorem~\ref{thm:subspace_bound}. % We adapt the approach of Henk,
                                % Sch\"urmann and
                                % Wills~\cite{Henk:2005}, who proved
                                % the assertion (without the equality
                                % condition) in the special case where
                                % $q=n$ and $K$ is a symmetric
                                % polytope, only using the
                                % Brunn-Minkowski theorem. 
We use Fubini's theorem to decompose the dual curvature measure into
integrals over hyperplane sections. 
Lemma~\ref{lem:unimodal} will provide a critical estimate for these integrals.
\begin{proof}[Proof of Theorem~\ref{thm:subspace_bound}] 
In order to prove the inequality \eqref{eq:subspace_bound} wer
we may certainly assume  $q>\dim L=k$. For $\vy\in K|L$ put $\ov{\vy}=\rho_{K|L}(\vy)\vy$ and
\begin{align*}
	F_\vy &= \conv\{\vnull, K\cap(\ov{\vy}+L^\perp)\},\\
	M_\vy &= \conv\{K\cap L^\perp , K\cap(\ov{\vy}+L^\perp)\}.
\end{align*}
Observe that $	M_\vy\cap (\ov{\vy}+L^\perp)= F_\vy\cap
(\ov{\vy}+L^\perp) = K\cap(\ov{\vy}+L^\perp)$. 
By Lemma~\ref{lem:curvature_measure_euclidean_coordinates} and
Fubini's theorem we may write 
\begin{equation}
\begin{split} 
	\dcura_q(K,S^{n-1}) &= \frac{q}{n} \int\limits_{K|L} \left(\,\int\limits_{K\cap (\vy+L^\perp)} |\vz|^{q-n} \mathrm{d}\mathcal{H}^{n-k}(\vz) \right)\mathrm{d}\mathcal{H}^k(\vy)\\
	&\geq \frac{q}{n} \int\limits_{K|L} \left(\,
          \int\limits_{M_\vy\cap (\vy+L^\perp)} |\vz|^{q-n}
          \mathrm{d}\mathcal{H}^{n-k}(\vz)
        \right)\mathrm{d}\mathcal{H}^k(\vy).
\end{split}
\label{eq:first_step} 
\end{equation}
In order to estimate the inner integral let  $\vy\in K|L$, $\vy\neq \vnull$, and for abbreviation we set
$\lambda=\rho_{K|L}(\vy)^{-1}\leq 1$. Then by the symmetry of $K$ we find 
\begin{align*}
	M_\vy\cap(\vy+L^\perp) \supseteq& \lambda (K\cap(\ov{\vy}+L^\perp))+(1-\lambda)(K\cap L^\perp)\\
	\supseteq& \lambda (K\cap(\ov{\vy}+L^\perp))\\
	&+(1-\lambda)\left( \frac12 ( K\cap(\ov{\vy}+L^\perp) ) + \frac12 ( - ( K\cap(\ov{\vy}+L^\perp) ) ) \right)\\
	=& \frac{1+\lambda}{2}(K\cap
           (\ov{\vy}+L^\perp))+\frac{1-\lambda}{2}(-(K\cap
           (\ov{\vy}+L^\perp))). 
\end{align*}
Hence the  set $M_\vy\cap(\vy+L^\perp)$ contains a convex combination
of a set and its reflection at the origin. This allows us to apply Lemma~\ref{lem:unimodal} from which we get
\begin{equation}
	\int\limits_{M_\vy\cap(\vy+L^\perp)} |\vz|^{q-n} \mathrm{d}\mathcal{H}^{n-k}(\vz) \geq \int\limits_{K\cap (\ov{\vy}+L^\perp)} |\vz|^{q-n} \mathrm{d}\mathcal{H}^{n-k}(\vz).
	\label{eq:inner_integral_estimate}
\end{equation}
Together with \eqref{eq:first_step} we obtain the lower bound  
\begin{equation}
	\dcura_q(K,S^{n-1}) \geq \frac{q}{n} \int\limits_{K|L}\left(\, \int\limits_{K\cap (\ov{\vy}+L^\perp)} |\vz|^{q-n} \mathrm{d}\mathcal{H}^{n-k}(\vz)\right) \mathrm{d}\mathcal{H}^k(\vy).% \\
	% &= \frac{q}{n} \int\limits_{K|L}\left(\,
        % \int\limits_{F_\vy\cap (\ov{\vy}+L^\perp)} |\vz|^{q-n}
        % \mathrm{d}\mathcal{H}^{n-k}(\vz)\right)
        % \mathrm{d}\mathcal{H}^k(\vy).
\label{eq:lower_bound} 
\end{equation}	
On the other hand we find
\begin{align*}
	\dcura_q(K,S^{n-1}\cap L)&=\frac{q}{n} \int\limits_{K|L} \left(\,\int\limits_{F_\vy\cap(\vy+L^\perp)} |\vz|^{q-n} \mathrm{d}\mathcal{H}^{n-k}(\vz)\right)\mathrm{d}\mathcal{H}^k(\vy)\\
	&=\frac{q}{n} \int\limits_{K|L} \left(\,\int\limits_{(\rho_{K|L}(\vy)^{-1}(K\cap(\ov{\vy}+L^\perp))} |\vz|^{q-n} \mathrm{d}\mathcal{H}^{n-k}(\vz)\right)\mathrm{d}\mathcal{H}^k(\vy)\\
	&=\frac{q}{n} \int\limits_{K|L} \rho_{K|L}(\vy)^{k-q}\left( \int\limits_{K\cap(\ov{\vy}+L^\perp)} |\vz|^{q-n} \mathrm{d}\mathcal{H}^{n-k}(\vz)\right)\mathrm{d}\mathcal{H}^k(\vy).
\end{align*}
% where in going from the second to the third line we scaled $z$ by $\rho_{K|L}(y)^{-1}$.

The inner integral is independent of the length of $\vy\in K|L$ and might be as well considered as the value $g(\vu)$ of a (measurable) function $g\colon S^{n-1}\cap L\to\R_{\geq 0}$. By taking this into account and using spherical coordinates we obtain
%\begin{align*}
%	\frac{\dcura_q(K,S^{n-1}\cap L)}{\dcura_q(K,S^{n-1})} \leq& \frac{\int\limits_{K|L} \rho_{K|L}(y)^{k-q} g(y/|y|) \mathrm{d}\mathcal{H}^k(y)}{\int\limits_{K|L} g(y/|y|) \mathrm{d}\mathcal{H}^k(y)}\\
%	=& \frac{\int\limits_{S^{n-1}\cap L} g(u) \int\limits_0^{\rho_{K|L}(u)}\rho_{K|L}(ru)^{k-q} r^{k-1} \mathrm{d}r \mathrm{d}\mathcal{H}^{k-1}(u)}{\int\limits_{S^{n-1}\cap L} g(u) \int\limits_0^{\rho_{K|L}(u)} r^{k-1} \mathrm{d}r \mathrm{d}\mathcal{H}^{k-1}(u)}\\
%	=& \frac{\int\limits_{S^{n-1}\cap L} g(u) \cdot \frac1q\rho_{K|L}(u)^k \mathrm{d}\mathcal{H}^{k-1}(u)}{\int\limits_{S^{n-1}\cap L} g(u) \cdot \frac1k \rho_{K|L}(u)^k \mathrm{d}\mathcal{H}^{k-1}(u)}\\
%	&=\frac{k}{q}.
%\end{align*}
\begin{equation}
\begin{split} 
	\dcura_q(K,S^{n-1}\cap L) =& \frac{q}{n}\int\limits_{K|L} \rho_{K|L}(\vy)^{k-q} g(\vy/|\vy|) \mathrm{d}\mathcal{H}^k(\vy)\\
	=&\frac{q}{n}\int\limits_{S^{n-1}\cap L} g(\vu)
           \left(\int\limits_0^{\rho_{K|L}(\vu)}\rho_{K|L}(r\vu)^{k-q}
           r^{k-1} \mathrm{d}r\right)
           \mathrm{d}\mathcal{H}^{k-1}(\vu)\\
=&\frac{q}{n}\int\limits_{S^{n-1}\cap L} g(\vu) \rho_{K|L}(\vu)^{k-q}
           \left(\int\limits_0^{\rho_{K|L}(\vu)}r^{q-1} \mathrm{d}r\right)
           \mathrm{d}\mathcal{H}^{k-1}(\vu)\\
	=&\frac{1}{n}\int\limits_{S^{n-1}\cap L} g(\vu)
        \,\rho_{K|L}(\vu)^k \mathrm{d}\mathcal{H}^{k-1}(\vu).
\label{eq:numerator} 
\end{split} 
\end{equation}
Applying the same transformation to the right hand side of
\eqref{eq:lower_bound} gives 
\begin{equation}
\begin{split} 
	\dcura_q(K,S^{n-1}) \geq& \frac{q}{n}\int\limits_{K|L} g(\vy/|\vy|) \mathrm{d}\mathcal{H}^k(\vy)\\
	=& \frac{q}{n}\int\limits_{S^{n-1}\cap L} g(\vu) \left(\,\int\limits_0^{\rho_{K|L}(\vu)} r^{k-1} \mathrm{d}r\right) \mathrm{d}\mathcal{H}^{k-1}(\vu)\\
	=& \frac{q}{n}\frac{1}{k}\int\limits_{S^{n-1}\cap L} g(\vu)
        \,\rho_{K|L}(\vu)^k \mathrm{d}\mathcal{H}^{k-1}(\vu).
\label{eq:denominator}
\end{split}
\end{equation}
Combining \eqref{eq:numerator} and \eqref{eq:denominator} yields
\eqref{eq:subspace_bound} in the case $k=\dim L\leq q$, i.e., 
\begin{equation}
\frac{	\dcura_q(K,S^{n-1}\cap
  L)}{\dcura_q(K,S^{n-1})} \leq \frac{k}{q}.
\end{equation}

Now suppose that the dual curvature measure of $K$ satisfies
the inequality \eqref{eq:subspace_bound} 
with equality with respect to a proper subspace $L$. 
Then we certainly have $\dim L< q$, since the curvature measure cannot
be concentrated on a great hypersphere. Hence we must have equality in
\eqref{eq:inner_integral_estimate} for every $\vy\in\relint
(K|L)$. 

Assume $q< n$. Then the  superlevel sets $L_f^+(\alpha)$ of the function
$f(\vz)=|\vz|^{q-n}$, $\vz\in\R^n$, are balls. Hence, in view of the  equality condition of
Lemma~\ref{lem:unimodal}, equality in \eqref{eq:inner_integral_estimate} implies that  
$M_\vy\cap (\vy+L^\perp)\cap r B_n$ and $K\cap (\ov{\vy}+L^\perp)\cap r
B_n$ have the same $(n-k)$-dimensional volume for every $r>0$. For
sufficiently small $r$, however, the intersection of $K\cap
(\ov{\vy}+L^\perp)\subset \partial K$ with $r\,B_n$ is empty. Hence we must
have $q=n$ and in this case we know by Theorem
\ref{thm:logarithmic_minkowski_problem}  that equality 
% Let $x\in M_y\cap (\ov{y}+L^\perp)$ be a point of minimal length and put $r=|x|>0$. Then $\rho_{K|L}(y)^{-1}x$ lies in $M_y\cap (y+L^\perp)\cap\inte(rB_n)$, but $M_y\cap (\ov{y}+L^\perp)\cap\inte(rB_n)$ is empty.
% When $q=n$, equality in the subspace bound
is attained if and only if the cone-volume measure of $K$ satisfies the subspace concentration condition as stated in Theorem~\ref{thm:logarithmic_minkowski_problem}.
\end{proof}

\begin{remark} 
It is worth noting, that the proof of Theorem~\ref{thm:subspace_bound}
only relies on the symmetry  of the function
$|\cdot|^{q-n}=\rho_{B_n}(\cdot)^{n-q}$, its homogeneity 
and the convexity of its unit ball. In fact, the ball  $B_n$ can be
replaced by
any symmetric convex body 
$M\in\Ken$ in the sense that 
\begin{equation}
\begin{split}
	\int\limits_{\alpha_K^\ast(S^{n-1}\cap L)} &\rho_M(\vu)^{n-q}
        \rho_K(\vu)^q \mathrm{d}\mathcal{H}^{n-1}(\vu) \leq\\ 
    &\frac{\dim L}{q} \int\limits_{S^{n-1}} \rho_M(\vu)^{n-q} \rho_K(\vu)^q \mathrm{d}\mathcal{H}^{n-1}(\vu),
\end{split} 
\label{eq:general} 
\end{equation}
where $L\subseteq\R^n$ is a subspace with $\dim L\leq q$. Observe, in this more
general setting, Lemma \ref{lem:curvature_measure_euclidean_coordinates} becomes 
\begin{equation}
\begin{split} 
\int\limits_{\alpha_K^\ast(\eta)} &\rho_M(\vu)^{n-q} \rho_K(\vu)^q
\mathrm{d}\mathcal{H}^{n-1}(\vu) \\  
=& \int\limits_{\alpha_K^\ast(\eta)}\rho_M(\vu)^{n-q}\left(\int\limits_0^{\rho_K(\vu)}
q\,r^{n-1\,}r^{q-n} \mathrm{d}r\right)
\mathrm{d}\mathcal{H}^{n-1}(\vu) \\ = &q\,\int\limits_{\vx\in K,\,\,\vx/|\vx|\in\alpha_K^\ast(\eta)} \rho_M(\vx)^{n-q} \mathrm{d}\mathcal{H}^n(\vx), 
\end{split} 
\end{equation}
and \eqref{eq:general} can be proved along the same lines as
Theorem~\ref{thm:subspace_bound} with $\rho_{B_n}(\cdot)$ replaced by $\rho_M(\cdot)$.
\end{remark}

\section{Proof of  Proposition \ref{prop:bound_opt}}
Here we show that the bounds given in Theorem~\ref{thm:subspace_bound}
are indeed tight for every choice of  $q\in(0,n)$. To this end let
$k\in\N$ with $0<k<n$ and for $r>0$ let $K_r$ be the cylinder 
\begin{equation}
  K_r=(rB_k)\times B_{n-k}. 
\end{equation} 
Let  $L=\lin \{ \ve_1,\ldots,\ve_k \}$ be the $k$-dimensional subspace
generated by  the first $k$ canonical unit vectors $\ve_i$. 

%\begin{example}
% Let $k,n\in\N$, $k<n$, and $q\in(0,n]$. For $r>0$ put $K=(rB_k)\times B_{n-k}$ and $L=\lin \{ e_1,\ldots,e_k \}$, where the vectors $e_1,\ldots,e_n$ denote the canonical basis of $\R^n$.
For $\vx\in \R^n$ write $\vx=\vx_1+\vx_2$, where
$\vx_1\in\R^k\times\{\vnull\}$ and
$\vx_2\in\{\vnull\}\times\R^{n-k}$. The supporting hyperplane of $K_r$
with respect to  a unit vector $\vv\in S^{n-1}\cap L$ is given by
\begin{equation}
	H_{K_r}(\vv)=\{ \vx\in \R^n\colon \langle \vv_1,\vx_1\rangle=r \}.
\end{equation}
Hence the part of the boundary of $K_r$ covered by all these supporting
hyperlanes is given by $rS^{k-1}\times B_{n-k}$. % and 
% This shows that $\alpha_K^\ast(S^{n-1}\cap L)$ consists of all points $x$ with $r|x_2|\leq |x_1|\leq r$.
In view of Lemma~\ref{lem:curvature_measure_euclidean_coordinates} and
Fubini's theorem we conclude 
%\begin{equation}
%	\frac{\dcura_q(K,S^{n-1}\cap L)}{\dcura_q(K,S^{n-1})}=\frac{\int\limits_{x_1\in rB_k} \int\limits_{\substack{x_2\in B_{n-k}\\r|x_2|\leq |x_1|}} (|x_1|^2+|x_2|^2)^{\frac{q-n}{2}} \mathrm{d}\mathcal{H}^{n-k}(x_2) \mathrm{d}\mathcal{H}^{k}(x_1)}{\int\limits_{x_1\in rB_k} \int\limits_{x_2\in B_{n-k}} (|x_1|^2+|x_2|^2)^{\frac{q-n}{2}} \mathrm{d}\mathcal{H}^{n-k}(x_2) \mathrm{d}\mathcal{H}^{k}(x_1)}.
%\end{equation}
\begin{equation}
\begin{split} 
	\dcura_q(K_r,\,&S^{n-1}\cap L) =\\  &\frac{q}{n} \int\limits_{\vx_1\in rB_k} \left(\,\int\limits_{\substack{\vx_2\in B_{n-k}\\r|\vx_2|\leq |\vx_1|}} (|\vx_1|^2+|\vx_2|^2)^{\frac{q-n}{2}} \mathrm{d}\mathcal{H}^{n-k}(\vx_2) \right)\mathrm{d}\mathcal{H}^{k}(\vx_1).
	\label{eq:curvature_measure_cylinder_fubini}
\end{split}
\end{equation}
Denote the volume of $B_n$ by $\omega_n$. Recall, that the surface
area of $B_n$ is given by $n\omega_n$ and for abbreviation we set  
\begin{equation} 
c=c(q,k,n)=\frac{q}{n} k\omega_k (n-k)\omega_{n-k}.
\end{equation}
Switching to the cylindrical coordinates 
\begin{equation}
		\vx_1=s\vu,\quad s\geq 0, \vu\in S^{k-1},\qquad 
		\vx_2=t\vv,\quad t\geq 0, \vv\in S^{n-k-1},
	% \mathrm{d}\mathcal{H}^{n-k}(\vx_2) \mathrm{d}\mathcal{H}^{k}(\vx_1) = s^{i-1} t^{n-i-1} \mathrm{d}t \mathrm{d}s \mathrm{d}\mathcal{H}^{k-1}(u) \mathrm{d}\mathcal{H}^{n-k-1}(v)
\end{equation}
% \begin{gather*}
% 	\begin{aligned}[t]
% 		\vx_1&=s\vu,\qquad&&s\geq 0, \vu\in S^{k-1},\\
% 		\vx_2&=t\vv,\qquad&&t\geq 0, \vv\in S^{n-k-1},
% 	\end{aligned}\\
% 	% \mathrm{d}\mathcal{H}^{n-k}(\vx_2) \mathrm{d}\mathcal{H}^{k}(\vx_1) = s^{i-1} t^{n-i-1} \mathrm{d}t \mathrm{d}s \mathrm{d}\mathcal{H}^{k-1}(u) \mathrm{d}\mathcal{H}^{n-k-1}(v)
% \end{gather*}
transforms the right hand side of
\eqref{eq:curvature_measure_cylinder_fubini} to 
% the cylindrical coordinate system
% \begin{gather*}
% 	\begin{aligned}[t]
% 		x_1&=su,\qquad&&s\geq 0, u\in S^{k-1},\\
% 		x_2&=tv,\qquad&&t\geq 0, v\in S^{n-k-1},
% 	\end{aligned}\\
% 	\mathrm{d}\mathcal{H}^{n-k}(x_2) \mathrm{d}\mathcal{H}^{k}(x_1) = s^{i-1} t^{n-i-1} \mathrm{d}t \mathrm{d}s \mathrm{d}\mathcal{H}^{k-1}(u) \mathrm{d}\mathcal{H}^{n-k-1}(v)
% \end{gather*}
% and appropriate rescaling of the variables yields
%\begin{align*}
%	\frac{\dcura_q(K,S^{n-1}\cap L)}{\dcura_q(K,S^{n-1})}&=\frac{\int\limits_0^r \int\limits_0^{s/r} s^{k-1} t^{n-k-1} (s^2+t^2)^{\frac{q-n}{2}} \mathrm{d}t \mathrm{d}s}{\int\limits_0^r \int\limits_0^1 s^{k-1} t^{n-k-1} (s^2+t^2)^{\frac{q-n}{2}} \mathrm{d}t \mathrm{d}s}\\
%	&=\frac{r^{k-n} \int\limits_0^r \int\limits_0^1 s^{q-1} t^{n-k-1} (1+r^{-2}t^2)^{\frac{q-n}{2}} \mathrm{d}t \mathrm{d}s}{\int\limits_0^r \int\limits_0^1 s^{k-1} t^{n-k-1} (s^2+t^2)^{\frac{q-n}{2}} \mathrm{d}t \mathrm{d}s}\\
%	&=\frac{\int\limits_0^1 \int\limits_0^1 s^{q-1} t^{n-k-1} (r^2+t^2)^{\frac{q-n}{2}} \mathrm{d}t \mathrm{d}s}{\int\limits_0^1 \int\limits_0^1 s^{k-1} t^{n-k-1} (r^2s^2+t^2)^{\frac{q-n}{2}} \mathrm{d}t \mathrm{d}s}.
%\end{align*}
\begin{align}
	\notag
	\dcura_q(K_r,S^{n-1}\cap L) &= c \int\limits_0^r \int\limits_0^{s/r} s^{k-1} t^{n-k-1} (s^2+t^2)^{\frac{q-n}{2}}\, \mathrm{d}t\, \mathrm{d}s\\
	\notag
	&= c \int\limits_0^r \int\limits_0^1 s^{q-1} r^{k-n} t^{n-k-1} (1+r^{-2}t^2)^{\frac{q-n}{2}} \mathrm{d}t\, \mathrm{d}s\\
	&= c\, r^k \int\limits_0^1 \int\limits_0^1 s^{q-1} t^{n-k-1} (r^2+t^2)^{\frac{q-n}{2}} \mathrm{d}t\, \mathrm{d}s.
	\label{eq:curvature_measure_cylinder_L}
\end{align}
Analogously we obtain
\begin{equation}
\begin{split} 
	\dcura_q(K_r,S^{n-1}) 
& = \frac{q}{n} \int\limits_{\vx_1\in rB_k}
\left(\,\int\limits_{\vx_2\in B_{n-k}}
  (|\vx_1|^2+|\vx_2|^2)^{\frac{q-n}{2}}
  \mathrm{d}\mathcal{H}^{n-k}(\vx_2)
\right)\mathrm{d}\mathcal{H}^{k}(\vx_1) \\
&  = c \int\limits_0^r \int\limits_0^{1} s^{k-1} t^{n-k-1} (s^2+t^2)^{\frac{q-n}{2}}\, \mathrm{d}t\, \mathrm{d}s\\
& = c\, r^k \int\limits_0^1 \int\limits_0^1 s^{k-1} t^{n-k-1} (r^2s^2+t^2)^{\frac{q-n}{2}} \mathrm{d}t\, \mathrm{d}s.
	\label{eq:curvature_measure_cylinder}
\end{split} 
\end{equation}
When $q>k$, the monotone convergence theorem gives
%\begin{align*}
%	\lim_{r\to 0+}\frac{\dcura_q(K,S^{n-1}\cap L)}{\dcura_q(K,S^{n-1})} &= \lim_{r\to 0+} \frac{\int\limits_0^1 \int\limits_0^1 s^{q-1} t^{n-k-1} (r^2+t^2)^{\frac{q-n}{2}} \mathrm{d}t \mathrm{d}s}{\int\limits_0^1 \int\limits_0^1 s^{k-1} t^{n-k-1} (r^2s^2+t^2)^{\frac{q-n}{2}} \mathrm{d}t \mathrm{d}s}\\
%	&= \frac{\int\limits_0^1 s^{q-1} \mathrm{d}s \cdot \int\limits_0^1 t^{q-k-1} \mathrm{d}t}{\int\limits_0^1 s^{k-1} \mathrm{d}s \cdot \int\limits_0^1 t^{q-k-1} \mathrm{d}t}=\frac{k}{q}.
%\end{align*}
\begin{equation}
	\lim_{r\to 0+} \int\limits_0^1 \int\limits_0^1 s^{q-1} t^{n-k-1} (r^2+t^2)^{\frac{q-n}{2}} \mathrm{d}t\, \mathrm{d}s = \int\limits_0^1 s^{q-1} \mathrm{d}s \cdot \int\limits_0^1 t^{q-k-1} \mathrm{d}t = \frac{1}{q(q-k)}
\end{equation}
and
\begin{equation}
\lim_{r\to 0+} \int\limits_0^1 \int\limits_0^1 s^{k-1} t^{n-k-1} (r^2s^2+t^2)^{\frac{q-n}{2}} \mathrm{d}t\,\mathrm{d}s = \int\limits_0^1 s^{k-1} \mathrm{d}s \cdot \int\limits_0^1 t^{q-k-1} \mathrm{d}t = \frac{1}{k(q-k)}.
\end{equation}
Hence, by \eqref{eq:curvature_measure_cylinder_L} and
\eqref{eq:curvature_measure_cylinder} we get 
\begin{equation}
	\lim_{r\to 0+}\frac{\dcura_q(K_r,S^{n-1}\cap L)}{\dcura_q(K_r,S^{n-1})} = \frac{k}{q}.
\end{equation}
%Now suppose $q\leq k$. Then
%\begin{align*}
%\frac{\dcura_q(K,S^{n-1}\cap L)}{\dcura_q(K,S^{n-1})} &= \frac{\int\limits_0^1 \int\limits_0^1 s^{q-1} t^{n-k-1} (r^2+t^2)^{\frac{q-n}{2}} \mathrm{d}t \mathrm{d}s}{\int\limits_0^1 \int\limits_0^1 s^{k-1} t^{n-k-1} (r^2s^2+t^2)^{\frac{q-n}{2}} \mathrm{d}t \mathrm{d}s}\\
%&= \frac{\int\limits_0^1 \int\limits_0^1 s^{q-1} t^{n-k-1} (r^2+t^2)^{\frac{q-n}{2}} \mathrm{d}t \mathrm{d}s}{\int\limits_0^1 \int\limits_0^{1/s} s^{k-1} s^{n-k-1} t^{n-k-1} (r^2s^2+s^2t^2)^{\frac{q-n}{2}} s \mathrm{d}t \mathrm{d}s}\\
%&= \frac{\int\limits_0^1 \int\limits_0^1 s^{q-1} t^{n-k-1} (r^2+t^2)^{\frac{q-n}{2}} \mathrm{d}t \mathrm{d}s}{\int\limits_0^1 \int\limits_0^{1/s} s^{q-1} t^{n-k-1} (r^2+t^2)^{\frac{q-n}{2}} \mathrm{d}t \mathrm{d}s}\\
%&= 1- \frac{\int\limits_0^1 \int\limits_1^{1/s} s^{q-1} t^{n-k-1} (r^2+t^2)^{\frac{q-n}{2}} \mathrm{d}t \mathrm{d}s}{\int\limits_0^1 \int\limits_0^{1/s} s^{q-1} t^{n-k-1} (r^2+t^2)^{\frac{q-n}{2}} \mathrm{d}t \mathrm{d}s}.
%\end{align*}
Now suppose $q\leq k$. Rewrite \eqref{eq:curvature_measure_cylinder} as
\begin{equation}
	\dcura_q(K_r,S^{n-1}) = c\,r^k \int\limits_0^1
        \int\limits_0^{1/s} s^{q-1} t^{n-k-1}
        (r^2+t^2)^{\frac{q-n}{2}} \mathrm{d}t\, \mathrm{d}s,
\label{eq:rewrite} 
\end{equation}
which in view of  \eqref{eq:curvature_measure_cylinder_L} gives
\begin{equation}
\begin{split} 
\dcura_q(K_r,S^{n-1}\cap L) & -
\dcura_q(K_r,S^{n-1})\\  
&=c\,r^k \int\limits_0^1\int\limits_1^{1/s} s^{q-1} t^{n-k-1} (r^2+t^2)^{\frac{q-n}{2}}
\mathrm{d}t\, \mathrm{d}s.
\end{split} 
\label{eq:diff} 
\end{equation}
Observe, that by the monotone convergence theorem
\begin{equation}
\begin{split} 
	\lim_{r\to 0+} \int\limits_0^1 \int\limits_1^{1/s} s^{q-1} t^{n-k-1} (r^2+t^2)^{\frac{q-n}{2}} \mathrm{d}t\, \mathrm{d}s &= \int\limits_0^1 \int\limits_1^{1/s} s^{q-1} t^{q-k-1} \mathrm{d}t\, \mathrm{d}s\\
	&= \begin{cases}
		\int_0^1 s^{q-1} \frac{1-s^{k-q}}{k-q} \mathrm{d}s, & \text{ if $q<k$,}\\
		\int_0^1 s^{q-1}(-\log s) \mathrm{d}s, & \text{ if $q=k$,}
	\end{cases}\\
	&= \frac{1}{kq}.
\label{eq:bounded}
\end{split}
\end{equation}
On the other hand, if $0<r<1$, then
\begin{equation}
\begin{split}
	\int\limits_0^1 \int\limits_0^{1/s} s^{q-1} t^{n-k-1} (r^2+t^2)^{\frac{q-n}{2}} \mathrm{d}t\, \mathrm{d}s &\geq \int\limits_0^1 \int\limits_r^1 s^{q-1} t^{n-k-1} (r+t)^{q-n} \mathrm{d}t \mathrm{d}s\\
	&\geq \int\limits_0^1 \int\limits_r^1 s^{q-1} t^{n-k-1} (t+t)^{q-n} \mathrm{d}t\, \mathrm{d}s\\
	&= \begin{cases}
		\frac{2^{q-n}}{q} \frac{r^{q-k}-1}{k-q}, & \text{if $q<k$,}\\
		\frac{2^{q-n}}{q} (-\log(r)), & \text{if $q=k$,}
	\end{cases}
\label{eq:unbounded} 
\end{split}
\end{equation} 
which is not bounded from above as a function in $r$. Hence, by
\eqref{eq:rewrite}, \eqref{eq:diff}, \eqref{eq:bounded}, \eqref{eq:unbounded} we finally get 
\begin{equation}
	\lim_{r\to 0+}\frac{\dcura_q(K_r,S^{n-1}\cap
          L)}{\dcura_q(K_r,S^{n-1})} = 1- \lim_{r\to 0+} \frac{\int\limits_0^1\int\limits_1^{1/s} s^{q-1} t^{n-k-1} (r^2+t^2)^{\frac{q-n}{2}}
\mathrm{d}t\, \mathrm{d}s}{	\int\limits_0^1 \int\limits_0^{1/s} s^{q-1} t^{n-k-1} (r^2+t^2)^{\frac{q-n}{2}} \mathrm{d}t\, \mathrm{d}s}=1,
\end{equation}
which finishes the proof of Proposition \ref{prop:bound_opt}.
%\end{example}

%\bigskip
%\noindent 
%{\it  Acknowledgement.} The authors thank  Rolf Schneider,
%Eugenia Saor{\'\i}n G{\'o}mez, Guangxian Zhu and the referees for their very helpful
%comments and suggestions.
%\bigskip
%\bigskip

%\bibliographystyle{plain}
%\bibliography{references}

\end{document}